\documentclass[12pt]{article}
\usepackage{amsbsy,amsfonts}

\newtheorem{theo}{Theorem}[section]
\newtheorem{lemma}[theo]{Lemma}

\newtheorem{coll}[theo]{Corollary}
\newtheorem{remark}[theo]{Remark}

\newcommand{\be}{\begin{equation}}
\newcommand{\ee}{\end{equation}}
\newcommand{\ba}{\begin{array}}
\newcommand{\ea}{\end{array}}

\newcommand{\m}{\frak m}
\newcommand{\dsum}{\displaystyle \sum}

\begin{document}

\title{ Counting of paths and the multiplicity of determinantal rings}
\author{ Hsin-Ju Wang \thanks{e-mail:  hjwang@math.ccu.edu.tw} \\
Department of Mathematics, National Chung Cheng University\\
          Chiayi 621, Taiwan}

\date{}

\maketitle

\begin{abstract}
~\\ In this paper, we derive several formulas of counting families
of non-intersecting paths for two-sided ladder-shaped regions. As
an application we give a new proof to a combinatorial
interpretation of Fibonacci numbers obtained by G. Andrews in
1974.
\end{abstract}

\section{Introduction}
Let $X=(x_{ij})$ be  an $m\times n$ matrix of indeterminates over
a field $K$. Let $I_{r+1}$ be the ideal of $K[X]$ generated by the
set of all $(r+1)$-minors of $X$ and set $R_{r+1}=K[X]/I_{r+1}$;
then it is well known that the multiplicity of $R_{r+1}$ is given
by $$ e(R_{r+1})=\det[{m+n-i-j\choose m-i}]_{i, j=1, \dots, r}.$$
The formula was first found by G.Z. Giambelli \cite{gi} in 1909,
and a different formula was given by Abhyankar \cite{ab} and
Galligo \cite{Ga} later. Recently, Herzog and Trung generalized
the formula as follows. \begin{theo} Let $X=(x_{ij})$ be a generic
$m\times n$ matrix over a field $K$. Let $a_1< \cdots < a_r\leq
m$, $b_1< \cdots < b_r\leq n$ be positive integers. Let $D_t$
denote the part of $X$ consisting of the first $a_t-1$ rows and
the first $b_t-1$ columns. Let $G_t$ be the set of all $t$-minors
of $D_t$, $t=1, \dots, r$, and let $G_{r+1}$ be the set of all
$(r+1)$-minors of $X$. Let $I$ be the ideal of $K[X]$ generated by
$G=\cup_{t=1}^{r+1} G_t$; then the multiplicity of $K[X]/I$ is
$$\det[{m+n-a_i-b_j\choose m-a_i}]_{i, j=1, \dots, r}.$$
\label{theo1}
\end{theo}

The proof of Theorem~\ref{theo1} can be sketched as follows. It is
well known, from the theory of Gr\"{o}bner bases, that the Hilbert
function of a homogeneous ideal $I$ in a polynomial ring $K[X]$
coincides with the one of the ideal $I^{\ast}$ generated by the
leading terms of the polynomials of $I$ (with respect to some term
order). Therefore we can replace $I$ by $I^{\ast}$ which is a
monomial ideal. Since $I^{\ast}$ is a square free monomial ideal,
we can associate with $I^{\ast}$ a simplicial complex $\Delta$
such that $K[X]/I^{\ast}$ is $K[\Delta]$, the face ring of
$\Delta$. By a result of R. Stanley \cite{st}, the multiplicity of
$K[\Delta]$ is equal to the number of facets (faces of maximal
dimensions) of $\Delta$ which can be characterized as families of
non-intersecting paths. Therefore, the formula can be obtained by
a method of I. Gessel and G. Viennot \cite{gevi} (see also
\cite[Sect.~2.7]{st1}). \par  The goal of this paper is to compute
the multiplicity of certain ladder determinantal rings. Ladder
determinantal rings arise for instance  in Abhyankar's study on
the singularities of Schubert varities of flag manifold. Let's
recall the definition of ladder determinantal rings. Let
$X=(x_{ij})$ be an $m\times n$ matrix of indeterminates over a
field $K$. A ladder of $X$ is a subset $Y$ which satisfies the
condition that whenever $x_{ij}, x_{i'j'}\in Y$ with $i<i'$ and
$j<j'$, then $x_{ij'}, x_{i'j}\in Y$. Let $Y$ be a ladder of $X$.
Let $G$ be the set of $(r+1)$-minors of $Y$ and $I$ be the ideal
of $K[Y]$ generated by $G$. It is shown in
\cite[Corollary~4.2]{wa} that $G$ is a Gr\"{o}bner basis for $I$
with respect to some term order of $K[Y]$, therefore $I^{\ast}$ is
generated by the leading terms of the polynomials of $G$. As a
consequence, the multiplicity of $K[Y]/I$ is equal to the number
of facets of a simplicial complex $\Delta$  associated to
$I^{\ast}$. Moreover, we shall see in section 3 that every facet
of $\Delta$ can be characterized as a family of non-intersecting
paths of certain ladder-shaped region of the plane. Therefore, the
computation of the multiplicity for $K[Y]/I$ boils down to
counting families of non-intersecting paths. In section 2, we
prove a {\it beautiful} formula of counting families of
non-intersecting paths for two-sided ladder-shaped regions, we
state it in the following.
\begin{theo}
\label{theo2} Let $m, n$ and $r$ be positive integers with $r\leq
m\leq n$. Let $X=\{(i, j)~|~1\leq i\leq m, 1\leq j\leq n\}$. We
define a partial order on $X$ by setting $(i, j)\leq (i', j')$ if
$i\geq i'$ and $j\leq j'$. A path $C$ from $P$ to $Q$ is a maximal
chain in $X$ with endpoints $P$ and $Q$. We use $w(P, Q)$ for the
number of all different paths from $P$ to $Q$. Let $Y$ be a
two-sided sub-region of $X$ (see Figure 1 for an example of $Y$).
Let $P_i, Q_i$, $i=1, \dots, r$, be points of $Y$ with the
following properties:
\begin{description}
\item{(i)} $P_i=(a_i, n)$ and $1=a_1<a_2< \cdots < a_r\leq m$.
\item{(ii)} $Q_i=(c_i, d_i)$, $1\leq c_1\leq \cdots \leq c_r=m$
and $1=d_1< d_2< \cdots < d_r\leq n$. \item{(iii)} If $(1,
1)\notin Y$, then $(c_1-1, 1)\in Y$.\end{description} Then the
number of non-intersecting paths from $P_i$ to $Q_i$, $i=1, \dots,
r$, is $$ \det[w(P_i, Q_i)]_{i, j=1, \dots, r},$$ where $w(P_i,
Q_i)$ is the set of all paths from $P_i$ to $Q_i$.
\end{theo} With the help of this formula and
Theorem~\ref{theo7}, we are able to obtain the following results.

\begin{theo}
Let $m, n, l, k$ and $r$ be non-negative integers  with $k, l\leq
min\{m-2, n-2\}$ and $r< min\{m-k, n-l\}$. Let $K$ be a field and
$X=(x_{ij})$ be an $m\times n$ matrix of indeterminates. Let
$$Y=\{x_{ij}~|~ 1\leq i\leq m, max\{1, l-i+2\}\leq j\leq min\{n,
n+m-k-i\}\}$$ be a ladder of $X$. Let $I=I_{r+1}(Y)$ be the ideal
of $K[Y]$ generated by the $(r+1)$-minors of $Y$ and set
$R=K[Y]/I$. Then the multiplicity of $R$ is $$\det[\sum_{(a, b)\in
T} (-1)^{a+b} {m+n-i-j\choose m-1+(a-b-1)(i-1)+a(k-m)+b(l-n)}
]_{i, j=1, \dots, r}, $$ where $T=\{(a, b)\in \mathbb{Z}\times
\mathbb{Z} ~|~a-b=0~or~1\}$.  \label{theo3}
\end{theo}
\begin{theo}
Let $r, l$ and $n$ be positive integers with $r<l<n$ and $2r\leq
n$. Let $X$ be a skew-symmetric $n\times n$ matrix of
indeterminates $x_{ij}$, over a field $K$. Let $$Y=\{x_{ij}~|~
max\{1, i-l\}\leq j\leq min\{n, l+i\}, 1\leq i\leq n\}$$ be a
diagonal ladder of $X$. Let $Q_r$ be the ideal of $K[Y]$ generated
by the set of all $2r$-pfaffians of $Y$ and let $R=K[Y]/Q_r$; then
the multiplicity of $R$ is $$\det[\sum_{(a, b)\in T} (-1)^{a+b}
{2n-2r-i-j\choose n-r-i+a(-r+i-1)+b(r-l-i)} ]_{i, j=1, \dots,
r-1}, $$ where $T=\{(a, b)\in \mathbb{Z}\times \mathbb{Z}
~|~a-b=0~or~1\}$. \label{theo4}
\end{theo}

\section{Formulas of counting families of non-intersecting paths}
Let $m, n$ and $r$ be positive integers with $r\leq m\leq n$.
Consider the set of points $X=\{(i, j)~|~1\leq i\leq m, 1\leq
j\leq n\}$ of the plane  We define a partial order on $X$ by
setting $(i, j)\leq (i', j')$ if $i\geq i'$ and $j\leq j'$. Let
$P, Q\in X$ with $P\geq Q$; a path $C$ from $P$ to $Q$ is a
maximal chain in $X$ with endpoints $P$ and $Q$. We use $w(P, Q)$
for the number of all different paths from $P$ to $Q$.\\ Let
$P_i$, $Q_i$, $i=1, \dots, r$, be points of $X$; a subset $W$ of
$X$ is called an $r$-tuple of non-intersecting  paths from $P_i$
to $Q_i$ ($i=1, \dots, r$) if $W=C_1\cup C_2\cup \cdots \cup C_r$
where each $C_i$ is a path from $P_i$ to $Q_i$ and where $C_i\cap
C_j=\emptyset$ if $i\neq j$. We use $w({\bf P}, {\bf Q})$ for the
number of all $r$-tuple of non-intersecting  paths from $P_i$ to
$Q_i$, where ${\bf P}=\{P_1, \dots, P_r\}$ and ${\bf Q}=\{Q_1,
\dots, Q_r\}$. If all the points of ${\bf P}$ are on the line
$y=n$ and all the points of ${\bf Q}$ are on the line $x=m$, then
we have the following result from \cite{st1}.

\begin{theo}
\cite[Sect~2.7]{st1} Let $X=\{(i, j)~|~1\leq i\leq m, 1\leq j\leq
n\}$ be a subset of the plane with the partial order given in the
beginning. Let $1\leq a_1< \cdots < a_r\leq m$ and $1\leq b_1<
\cdots < b_r\leq n$ be strictly increasing sequences of positive
integers. Let $P_i=(a_i, n)$ and $Q_i=(m, b_i)$; then the number
of all $r$-tuple of non-intersecting  paths from ${\bf P}$  to
${\bf Q}$  is $$\det[w(P_i, Q_j)]_{i,j=1, \dots, r}.$$
\label{theo5}
\end{theo}
\par The goal of this section is to show that the same equation
holds for ladder-shaped  regions. \par Let $C_{k_1,k_2}$ be a path
from $(k_1, n)$ to $(m, k_2)$ with $1<k_1\leq m$ and $1<k_2\leq
n$. Let $\tilde{C}_{k_3,k_4}$ be a path from $(1, k_3)$ to $(k_4,
1)$ with $1\leq k_3<n$ and $1\leq k_4<m$. If $Y$ is the sub-region
of $X$ bounded by $x=1$, $x=m$, $y=1$, $y=n$, $C_{k_1,k_2}$ and
$\tilde{C}_{k_3,k_4}$, then we say $Y$ is a two-sided sub-ladder
of $X$ (determined by $C_{k_1,k_2}$ and $\tilde{C}_{k_3,k_4}$).
Let $P_i$, $Q_i$, $i=1, \dots, r$, be points of $Y$ with the
following properties: \begin{description} \item{(i)} $P_i=(a_i,
n)$ and $1=a_1< a_2< \cdots <a_r\leq m$. \item{(ii)} $Q_i=(c_i,
d_i)$, $1\leq c_1\leq \cdots \leq c_r=m$ and $1=d_1<d_2< \cdots <
d_r\leq n$. \item{(iii)} $c_1>k_2$ if $k_2\geq 2$.
\end{description}  An example of $Y$ with $k_1=m-2$, $k_2=n-2$,
$k_3=4$ and $k_4=7$ is displayed in Figure 1. \vskip 0.2in
\begin{picture}(0,0)
\put(120,-10){\line(0, -1){70}}\put(120,-80){\line(1,0){20}}
\put(140,-80){\line(0,-1){10}} \put(140,-90){\line(1,0){10}}
\put(150,-90){\line(0,-1){10}} \put(150,-100){\line(1,0){30}}
\put(180,-100){\line(0,-1){10}} \put(180,-110){\line(1,0){40}}
\put(220,-110){\line(0,1){80}}

\put(120,-10){\line(1,0){80}} \put(200,-10){\line(0,-1){10}}
\put(200,-20){\line(1,0){10}} \put(210,-20){\line(0,-1){10}}
\put(210,-30){\line(1,0){10}}

\put(117,-13){$\times$} \put(105,-7){$P_1$}

\put(137,-13){$\times$} \put(125,-7){$P_2$}
\put(177,-13){$\times$} \put(165,-7){$P_r$}

\put(197,-13){$\times$} \put(197,-7){(m-2,n)}

\put(217,-33){$\times$} \put(220,-28){(m,n-2)}

\put(217,-63){$\times$} \put(222,-58){$Q_r$}

\put(217,-113){$\times$} \put(217,-117){(m,1)}

\put(177,-113){$\times$} \put(157,-120){(7,1)}
\put(192,-113){$\times$} \put(185,-120){$Q_1$}

\put(202,-103){$\times$} \put(198,-93){$Q_2$}

\put(117,-83){$\times$} \put(107,-90){(1,4)}

\put(145,-140){Figure~1}

\end{picture}
\vskip 2.3in We will show in the following that the same
conclusion of  Theorem~\ref{theo5} holds with respect to $Y$. In
fact, our result is more general as we allow the points of ${\bf
Q}$ to be in a more general position. However, before doing so, we
need a couple of lemmas.

\begin{lemma}
Let $Y$ be the region described as above with $(m-1, k_2)\in
C_{k_1,k_2}$. Let $P=(a, n)$ and $Q=(m, d)$; then $w(P,
Q)=\sum_{t=d}^{k_2} w(P, Q_t)$, where $Q_t=(m-1, t)$. \label{lem1}
\end{lemma}
\begin{proof}
Let $S$ be the set of all paths from $P$ to $Q$; then $S$ is the
disjoint union of $S_t$ for $t=d, \dots, k_2$, where $$ S_t=\{C\in
S~|~(m, t)\in C, (m, t+1)\notin C\}$$ if $t<k_2$ and $$
S_{k_2}=\{C\in S~|~(m, k_2)\in C \}.$$ For every $t$, let
$Q_t=(m-1, t)$ and let $S_t'$ be the set of all paths from $P$ to
$Q_t$; then there is a one to one correspondence between $S_t$ and
$S_t'$ (if $C\in S_t$ then $C-\{(m, d), \dots, (m, t)\}\in S_t'$),
it follows that $w(P, Q)=\sum_{t=d}^{k_2} w(P, Q_t)$.
\end{proof}

\begin{lemma}
Let $Y$ be the region described as above with $k_4=m-1$. Assume
that $(k_4, h)\in \tilde{C}_{k_3,k_4}$ and $(k_4, h+1)\notin
\tilde{C}_{k_3,k_4}$ for some $h\geq 2$. If $Q_i=(m, d_i)$ with
$1=d_1<d_2<d_3\leq h+1$; then $$ w(P, Q_2)-w(P, Q_3)=(d_3-d_2)w(P,
Q)$$ and $$ w(P, Q_1)-w(P, Q_2)=(d_2-1)w(P, Q),$$ where $P$ is any
point on the line $y=n$ and $Q=(m-1, h)=(k_4, h)$. \label{lem2}
\end{lemma}
\begin{proof}
We may assume that $(m-1, k_2)\in C_{k_1,k_2}$. Then by
Lemma~\ref{lem1} $w(P, Q)=\sum_{t=d_i}^{k_2} w(P, Q_t')$ for every
$i$, where $Q_t'=(m-1, t)$. However, from the assumptions, we see
that $w(P, Q_t')=w(P, Q)$ for $1\leq t\leq h$, where $Q=(m-1, h)$,
therefore $$\ba{ccl} w(P, Q_2)-w(P, Q_3) & = & \sum_{t=d_2}^{k_2}
w(P, Q_t')- \sum_{t=d_3}^{k_2} w(P, Q_t') \\ & = &
\sum_{t=d_2}^{d_3-1} w(P, Q_t') \\ & = & (d_3-d_2)w(P, Q).\ea $$
Similarly, $w(P, Q_1)-w(P, Q_2)=(d_2-1)w(P, Q)$.
\end{proof}
\vskip 0.2in
\par Now we state the main result of this section as follows.
\begin{theo}
\label{theo6} Let $m, n$ and $r$ be positive integers with $r\leq
m\leq n$. Let $Y$ be the sub-region of $X$ determined by some
paths $C_{k_1,k_2}$ and $\tilde{C}_{k_3,k_4}$. Let $P_i$, $Q_i$,
$i=1, \dots, r$, be points of $Y$ with the following properties:
\begin{description} \item{(i)} $P_i=(a_i, n)$ and
$1=a_1<a_2< \cdots < a_r\leq m$. \item{(ii)} $Q_i=(c_i, d_i)$,
$1\leq c_1\leq \cdots \leq c_r=m$ and $1=d_1< d_2< \cdots <
d_r\leq n$. \item{(iii)} $c_1>k_2$ if $k_2\geq 2$.
\end{description} Then $$w({\bf P}, {\bf Q})=\det[w(P_i,
Q_j)]_{i,j=1, \dots, r}.$$ \end{theo}

\begin{remark}
\label{rem1} 1. Condition (ii) in the theorem says that no two
points of ${\bf Q}$ will be on the same horizontal line. \par 2.
If there is no path from $P_k$ to $Q_k$ for some $k$, then both
sides of the equations are $0$ since $w({\bf P}, {\bf Q})=0$ and
$w(P_i, Q_j)=0$ for $i\geq k$ and $j\leq k$.
\end{remark}
\vskip 0.5in {\bf Proof of Theorem 2.4}  Assume without loss of
generality that $(k_1, n-1), (m-1, k_2)\in C_{k_1,k_2}$ and $(2,
k_3), (k_4, 2)\in \tilde{C}_{k_3,k_4}$. We proceed the proof by
induction on $r$ and the area of the region $Y$. We shall consider
two situations: one is when $m-k_4\geq 2$ and the other is when
$m-k_4=1$. \par Assume that $m-k_4\geq 2$. In this case, we will
show the equation by induction on the number of the intersection
of ${\bf Q}$ with the vertical line $x=m$. So, let $l$ be the
non-negative integer with the property that $c_{l+1}=\cdots
=c_r=m$ but $c_l<m$, i.e., the intersection of ${\bf Q}$ with the
vertical line $x=m$ is $\{Q_{l+1}, \dots, Q_r\}$. \\ Assume that
$l=r-1$. If $P_r=(m, n)$ then $k_1=m$ and $k_2=n$. Furthermore,
$w(P_r, Q_r)=1$ and $w(P_r, Q_j)=0$ for $j=1, \dots, r-1$. Let
${\bf P'}=\{P_1, \dots, P_{r-1}\}$ and ${\bf Q'}=\{Q_1, \dots,
Q_{r-1}\}$; then by induction $w({\bf P'}, {\bf Q'})=\det[w(P_i,
Q_j)]_{i,j=1, \dots, r-1}$, it follows that $$w({\bf P}, {\bf Q})
= w({\bf P'}, {\bf Q'})= \det[w(P_i, Q_j)]_{i,j=1, \dots, r-1}=
\det[w(P_i, Q_j)]_{i,j=1, \dots, r}.$$ If $a_r<m$ and $Q_r=(m,
k_2)$, then $w(P_i, Q_r)=w(P_i, Q_r')$ for every $i$, where
$Q_r'=(m-1, k_2)$. Let $Q'_j=Q_j$ for $j=1, \dots, r-1$ and ${\bf
Q'}=\{Q_1', \dots, Q_r'\}$; then ${\bf P}$ and ${\bf Q'}$ are all
in a proper sub-region of $Y$, so that by induction $$ w({\bf P},
{\bf Q})= w({\bf P}, {\bf Q'})= \det[w(P_i, Q'_j)]_{i,j=1, \dots,
r}= \det[w(P_i, Q_j)]_{i,j=1, \dots, r}.$$ From the above, we may
assume that $a_r<m$ and $d_r<k_2$. Let $Q_{r,t}=(m, t)$ for
$t=d_r, \dots, k_2$. Let $S$ be the set of all $r$-tuple of
non-intersecting paths from $P_i$ to $Q_i$ ($i=1, \dots, r$); then
$S$ is the disjoint union of $S_t$ for $t=d_r, \dots, k_2$, where
$$S_t=\{C_1\cup C_2\cup \cdots \cup C_r\in S~|~Q_{r,t}\in C_r,
Q_{r,t+1}\notin C_r\}$$  \newpage   ~~\\ if $t\neq k_2$ and
$$S_{k_2}=\{C_1\cup C_2\cup \cdots \cup C_r\in S~|~Q_{r,k_2}\in
C_r\} . $$  For $t=d_r, \dots, k_2$, let $Q'_{r,t}=(m-1, t)$ and
let $S'_t$ be the set of all $r$-tuple of non-intersecting paths
from ${\bf P}$ to ${\bf Q'_t}$, where ${\bf Q'_t}=\{Q_1, \dots,
Q_{r-1}, Q_{r,t}\}$. It is clear that there is a one to one
correspondence between $S_t$ and $S_t'$ (if $C_1\cup C_2\cup
\cdots \cup C_r\in S_t$, then $C_1\cup C_2\cup \cdots \cup
C_{r-1}\cup (C_r-\{Q_{r,t},Q_{r,t-1}, \dots, Q_{r,d_r}\})\in
S'_t$), therefore $$|S_t|=|S_t'|=\det[w(P_i, Q_1) \cdots ~w(P_i,
Q_{r-1}) ~w(P_i, Q'_{r,t}]_{i=1, \dots, r}$$ by induction as ${\bf
P}$ and ${\bf Q'_t}$ are all in a proper sub-region of $Y$.
However, $\sum_{t=d_r}^{k_2} w(P_i, Q_{r,t}')=w(P_i, Q_r)$ by
Lemma~\ref{lem1}, it follows that $$w({\bf P}, {\bf
Q})=\dsum_{t=d_r}^{k_2} |S_t|= \det[w(P_i, Q_j)]_{i,j=1, \dots,
r}.$$ Assume that $l<r-1$. If $d_{l+2}-d_{l+1}=1$, then $w({\bf
P}, {\bf Q})=w({\bf P}, {\bf Q'})$, where ${\bf Q'}=\{Q_1, \dots,
Q_l, (m-1, d_{l+1}), Q_{l+2}, \dots, Q_r\}$. Let $Q_{l+1}'=(m-1,
d_{l+1})$; then it is clear that $w(P_i, Q_{l+1}')+w(P_i,
Q_{l+2})=w(P_i, Q_{l+1})$, it follows by induction that $$\ba{rcl}
& & w({\bf P}, {\bf Q})\\ & = & \det[w(P_i, Q_1) \cdots w(P_i,
Q_{l+1}')~w(P_i, Q_{l+2}) \cdots w(P_i, Q_r)]_{i=1, \dots, r}\\ &
= & \det[w(P_i, Q_1) \cdots w(P_i, Q_{l+1}')+w(P_i,
Q_{l+2})~w(P_i, Q_{l+2}) \cdots w(P_i, Q_r)]_{i=1, \dots, r}\\  &
= & \det[w(P_i, Q_1) \cdots w(P_i, Q_{l+1})~w(P_i, Q_{l+2}) \cdots
w(P_i, Q_r)]_{i=1, \dots, r}. \ea $$ Assume that
$d_{l+2}-d_{l+1}>1$. Let $Q_{l+1,t}=(m, t)$ for $t=d_{l+1}, \dots,
d_{l+2}-1$ ($Q_{l+1,d_{l+1}}=Q_{l+1}$). Let $S$ be the set of all
$r$-tuple of non-intersecting paths from $P_i$ to $Q_i$ ($i=1,
\dots, r$); then $S$ is the disjoint union of $S_t$ for
$t=d_{l+1}, \dots, d_{l+2}-1$, where $$S_t=\{C_1\cup \cdots C_r\in
S~|~ (m, t)\in C_{l+1}, (m, t+1)\notin C_{l+1}\}.$$

~~\\ For $t=d_{l+1}, \dots , d_{l+2}-1$, let $Q'_{l+1,t}=(m-1, t)$
and let $S'_t$ be the set of all $r$-tuple of non-intersecting
paths from ${\bf P}$ to ${\bf Q'_t}$, where $${\bf Q'_t}=\{Q_1,
\dots, Q_l, Q'_{l+1,t}, Q_{l+2}, \dots, Q_r\},$$ then as the above
$|S_t|=|S'_t|$ and $$|S'_t|=\det[w(P_i, Q_1) \cdots w(P_i,
Q'_{l+1,t})~w(P_i, Q_{l+2})\cdots w(P_i, Q_r)]_{i=1, \dots, r}$$
by induction. However, $$\ba{rcl} \sum_{t=d_{l+1}}^{d_{l+2}-1}
w(P_i, Q'_{l+1,t}) & = &  \sum_{t=d_{l+1}}^{k_2} w(P_i,
Q'_{l+1,t})- \sum_{t=d_{l+2}}^{k_2} w(P_i, Q'_{l+1,t}) \\ & =  &
w(P_i, Q_{l+1})-w(P_i, Q_{l+2}) \ea$$ by Lemma~\ref{lem1}, it
follows that $$\ba{rcl}  &  & w({\bf P}, {\bf Q}) \\ & = &
\sum_{t=d_{l+1}}^{d_{l+2}-1} |S_t| \\ & = &
\sum_{t=d_{l+1}}^{d_{l+2}-1} \det[w(P_i, Q_1) \cdots w(P_i,
Q'_{l+1,t})~w(P_i, Q_{l+2})\cdots w(P_i, Q_r)]_{i=1, \dots, r}\\ &
= & \det[w(P_i, Q_1) \cdots w(P_i, Q_{l+1})-w(P_i, Q_{l+2})~w(P_i,
Q_{l+2})\cdots w(P_i, Q_r) ]_{i=1, \dots, r}\\ & = & \det[w(P_i,
Q_j]_{i=1, \dots, r}. \ea $$ \par We now assume that $m-k_4=1$. If
this is the case, we have $Q_i=(m, d_i)$ and $Q_1=(m, 1)$.
Moreover, there exists an integer $h\geq 2$ such that $(k_4, h)\in
\tilde{C}_{k_3,k_4}$ and $(k_4, h+1)\notin \tilde{C}_{k_3,k_4}$.
To show the equation, we consider the following three situations:
(1). $d_3\leq h+1$, (2). $d_3>h+1$ and $d_2\leq h+1$, and (3).
$d_2>h+1$. \\ If $d_3\leq h+1$, then it is easy to see that
$w({\bf P}, {\bf Q})=0$. Moreover, by Lemma~\ref{lem2}, the first
three columns of the matrix $[w(P_i, Q_j)]_{i,j=1, \dots, r}$ are
linearly dependent, it follows that $\det[w(P_i, Q_j)]_{i,j=1,
\dots, r}=0= w({\bf P}, {\bf Q})$. \\ Suppose that $d_3> h+1$. Let
$S$ be the set of all $r$-tuple of non-intersecting paths from
$P_i$ to $Q_i$; then $S$ is the disjoint union of $S_t$ for $t=1,
\dots, d_2-1$, where $$S_t=\{C_1\cup \cdots \cup C_r\in S~|~ (m,
t)\in C_1, (m, t+1)\notin C_1\}.$$

~~\\ For $t=1, \dots, d_2-1$, let $Q'_{1,t}=(m-1, t)$ and let
$S'_t$ be the set of all $r$-tuple of non-intersecting paths from
${\bf P}$ to ${\bf Q'_t}$, where $${\bf Q'_t}=\{Q'_{1,t}, Q_2,
\cdots, Q_r\},$$ then as the above $|S_t|=|S'_t|$. Furthermore,
let $Q'_1=(m-1, h)$, $Q'_2=(m, h+1)$ and $S'$ be the set of all
$r$-tuple of non-intersecting paths from ${\bf P}$ to ${\bf Q'}$,
where $${\bf Q'}=\{Q'_1, Q'_2, Q_3, \cdots, Q_r\}.$$ Assume
$d_2\leq h+1$. Then it is easy to see that $|S'_t|=|S'|$ for every
$t$, hence $$w({\bf P}, {\bf Q})=\dsum_{t=1}^{d_2-1}
|S'_t|=(d_2-1)w({\bf P}, {\bf Q'}).$$ Since ${\bf P}$ and ${\bf
Q'}$ are all in a proper sub-region of $Y$ and the points of ${\bf
P}$ and ${\bf Q'}$ satisfy the conditions (i) to (iii), \be
\label{eq1} |S'|=w({\bf P}, {\bf Q'})=\det[w(P_i, Q'_1)~w(P_i,
Q'_2)~w(P_i, Q_3)\cdots w(P_i, Q_r)]_{i=1, \dots, r}\ee by
induction, therefore by Lemma~\ref{lem2} $$w(P, Q_1)-w(P_i,
Q_2)=(d_2-1)w(P_i, Q'_1)$$ and \be \label{eq2} w(P, Q_2)-w(P_i,
Q'_2)=(h+1-d_2)w(P_i, Q'_1),\ee we obtain that $$\ba{rcl}
 &  & \det[w(P_i, Q_1)~w(P_i,
Q_2)\cdots w(P_i, Q_r)]_{i=1, \dots, r}\\ & = & \det[w(P_i,
Q_1)-w(P_i, Q_2)~w(P_i, Q_2)\cdots w(P_i, Q_r)]_{i=1, \dots, r}\\
& = & \det[(d_2-1)w(P_i, Q'_1)~w(P_i, Q_2)\cdots w(P_i,
Q_r)]_{i=1, \dots, r}\\ & = & \det[(d_2-1)w(P_i, Q'_1)~w(P_i,
Q'_2)\cdots w(P_i, Q_r)]_{i=1, \dots, r} \\ & = & w({\bf P}, {\bf
Q}). \ea $$ Assume $d_2>h+1$. Then it is easy to see that
$|S'_t|=|S'|$ if $t\leq h$, hence \be \label{eq3} w({\bf P}, {\bf
Q})=\dsum_{t=1}^{d_2-1} |S'_t|=\dsum_{t=h+1}^{d_2-1}
|S'_t|+h|S'|.\ee

~~\\ Notice that ${\bf P}$ and ${\bf Q'_t}$ are all in a proper
sub-region of $Y$ for $t\geq h+1$, therefore, $$w({\bf P}, {\bf
Q'_t})= \det[w(P_i, Q'_{1,t})~w(P_i, Q_2) \cdots w(P_i,
Q_r)]_{i=1, \dots, r}$$ by induction. From (\ref{eq1}),
(\ref{eq2}), (\ref{eq3}) and the fact that $$ w(P_i, Q_1)-w(P_I,
Q_2)=hw(P_i, Q'_1)+\dsum_{h+1}^{d_2-1}w(P_i, Q'_{1,t}),$$ we
obtain that $$\ba{rcl} &  & \det[w(P_i, Q_1)~w(P_i, Q_2) \cdots
w(P_i, Q_r)]_{i=1, \dots, r}\\ & = & \det[w(P_i, Q_1)-w(P_i,
Q_2)~w(P_i, Q_2)\cdots w(P_i, Q_r)]_{i=1, \dots, r}  \\ & = &
\det[hw(P_i, Q'_1)~w(P_i, Q_2)\cdots w(P_i, Q_r)]_{i=1, \dots, r}
\\ &  & +\sum_{t=h+1}^{d_2-1} \det[w(P_i, Q'_{1,t})~w(P_i, Q_2)
\cdots w(P_i, Q_r)]_{i=1, \dots, r} \\ & = & h\det[w(P_i,
Q'_1)~w(P_i, Q'_2)~w(P_I, Q_3)\cdots w(P_i, Q_r)]_{i=1, \dots,
r}+\sum_{t=h+1}^{d_2-1}|S'_t| \\ & = &  w({\bf P}, {\bf Q}) .\ea
$$ This completes the proof. \vskip 0.2in Theorem\ref{theo6}
generalized the following result. \begin{coll}
\cite[Theorem~1]{gevi} Let $X=\{(i, j)~|~0\leq i\leq j\}$ be a
subset of the plane with the partial order given in the beginning.
Let $0\leq a_1< \cdots < a_k$ and $0\leq b_1< \cdots < b_k$ be
strictly increasing sequences of non-negative integers. Let
$A_i=(0, a_i)$ and $B_i=(b_i, b_i)$; then the number of $k$-tuple
of non-intersecting paths from $A_i$ to $B_i$, $i=1, \dots, k$, is
$$\det[{a_i\choose b_j}]_{i,j=1, \dots, k}.$$

\label{coll1}
\end{coll}
\begin{proof}
The assertion follows from the fact that $w(A_i, B_j)={a_i\choose
b_j}$.
\end{proof}

\section{Paths of two-sided regions}
Let $X=\{(i, j)~|~1\leq i\leq m, 1\leq j\leq n\}$. We define a
partial order on $X$ by setting ~$(i, j)\leq (i', j')$ ~if $i\geq
i'$ and $j\leq j'$. As before, if $P$ and $Q$ are points of $X$,
then a path $C$ from $P$ to $Q$ is a maximal chain in $X$ with end
points $P$ and $Q$. A path $C$ is called a {\it diagonal} path if
$C$ satisfies the following two conditions: (1). if $(i-1, j), (i,
j)\in C$ and $(i, j)$ is not the final point of $C$, then $(i,
j-1)\in C$. (2). If $(i, j+1), (i, j)\in C$ and $(i, j)$ is not
the final point of $C$, then $(i+1, j)\in C$. For convenience, let
$C_k$ denotes the diagonal path $\{(m-k+i, n-i)~|0\leq i\leq
k\}\cup \{(m-k+i, n-1-i~|~0\leq i\leq k-1\}$, and $\tilde{C}_l$
denotes the diagonal path $\{(i, l+2-i~|~1\leq i\leq l+1\}\cup
\{(i+1, l+2-i)~|~1\leq i\leq l\}$. \par Let $Y$ be a two-sided
sub-region of $X$. If $Y$ is determined by $C_k$ and $\tilde{C}_l$
for some integers $k$ and $l$, then we say that $Y$ is a two-sided
diagonal ladder (see Figure 2). The goal of this section is to
derive the following formula for such regions.

\begin{picture}(0,0)
\put(100,-10){\line(0, -1){40}}\put(100,-50){\line(1,0){10}}
\put(110,-50){\line(0,-1){10}}\put(110,-60){\line(1,0){10}}
\put(100,-10){\line(1,0){30}}\put(130,-10){\line(0,-1){10}}
\put(130,-20){\line(1,0){10}} \put(140,-20){\line(0,-1){10}}

\put(180,-60){\line(0,-1){30}}\put(180,-90){\line(-1,0){40}}
\put(140,-90){\line(0,1){10}}\put(140,-80){\line(-1,0){10}}
\put(130,-80){\line(0,1){10}}

\put(180,-60){\line(-1,0){10}}\put(170,-60){\line(0,1){10}}
\put(170,-50){\line(-1,0){10}}

\put(97,-13){$\times$} \put(82,-7){P(1,n)} \put(127,-13){$\times$}
\put(130,-7){(m-k,n)} \put(177,-63){$\times$}
\put(180,-60){(m,n-k)} \put(177,-78){$\times$}
\put(182,-72){$Q_t$} \put(167,-78){$\times$} \put(155,-72){$Q'_t$}

\put(177,-93){$\times$} \put(172,-100){Q(m,1)}
\put(137,-93){$\times$} \put(128,-100){(l+1,1)}
\put(97,-53){$\times$} \put(75,-62){(1,l+1)}

 \put(250,-80){Figure~2}
\end{picture}
\vskip 2.0in

\begin{theo}
Let $m, n, k$ and $l$ be non-negative integers with $k, l\leq
min\{m-2, n-2\}$. Let $X=\{(i, j)~|~1\leq i\leq m, 1\leq j\leq
n\}$ with the partial order given in the beginning. Let $Y$ be the
two-sided diagonal sub-ladder of $X$ determined by $C_k$ and
$\tilde{C}_l$. Let $P=(1, n)$ and $Q=(m, 1)$; then $$w(P, Q)=
\sum_{(i, j)\in T} (-1)^{i+j} {m+n-2\choose m-1+i(k-m)+j(l-n)}, $$
where $T=\{(i, j)\in \mathbb{Z}\times \mathbb{Z} ~|~i-j=0~or~1\}$.
\label{theo7}
\end{theo}
\begin{remark}
If $l=0$, i.e., $Y$ is an one-sided diagonal ladder of $X$
determined by $C_k$, then $$w(P, Q)={m+n-2\choose
m-1}-{m+n-2\choose k-1}$$  as $${m+n-2\choose m-1+i(k-m)-jn}=0$$
unless $(i, j)=(0, 0)$ or $(1, 0)$. \label{rem2}
\end{remark}
We need several lemmas before proving the formula. The first lemma
is fairly easy, we omit the proof.  \begin{lemma} Let $T=\{(i,
j)\in \mathbb{Z}\times \mathbb{Z} ~|~i-j=0~or~1\}$. Then the
following hold. \begin{description} \item{(i)} If $\phi_1: T
\longrightarrow T$ is a map given by $\phi_1((i, j))=(1-i, -j)$,
then $\phi_1$ is a bijection.
\item{(ii)} If $\phi_2: T
\longrightarrow T$ is a map given by $\phi_2((i, j))=(-i, -j-1)$,
then $\phi_2$ is a bijection.
\end{description}
\label{lem3}
\end{lemma}
\begin{lemma}
The following identities hold. \\ (i). $${m+k-2\choose
m-1+i(k-m)+j(l-n)}={m+k-2\choose m-1+(1-i)(k-m)-j(l-n)}.$$ (ii).
If $m=l+2$, then $$ {m+n-4\choose
m-3+i(k-m)+j(l-n)}={\m+n-4\choose m-3-i(k-m)-(j+1)(l-n)}.$$
\label{lem4}
\end{lemma}
\begin{proof}
The first equality follows from the fact that
$$[m-1+i(k-m)+j(l-n)]+[m-1+(1-i)(k-m)-j(l-n)]=m+k-2.$$ Similarly,
the second equality follows from the fact that $l=m-2$ and
$$[m-3+i(k-m)+j(l-n)]+[m-3-i(k-m)-(j+1)(l-n)]=m+n-4.$$
\end{proof}

\begin{lemma}
If $t>l$ and $k\leq n-t$, then $$\sum_{(i, j)\in T} (-1)^{i+j}
{m+n-t-2\choose m-2+i(k-m)+j(l-n)}={m+n-t-2\choose
m-2}-{m+n-t-2\choose k-2}.$$ \label{lem5}
\end{lemma}
\begin{proof}
To show the equation, it is enough to show that $${m+n-t-2\choose
m-2+i(k-m)+j(l-n)}=0$$ unless $j=0$. For this, let $(i, j)\in T$.
If $j<0$, then $i\leq 0$, so that $$m-2+i(k-m)+j(l-n)\geq
m+n-l-2>m+n-t-2,$$ it follows that $${m+n-t-2\choose
m-2+i(k-m)+j(l-n)}=0.$$ If $j> 0$, then $i\geq 1$, so that
$$m-2+i(k-m)+j(l-n)\leq m-2+(k-m)+(l-n)<k+t-n\leq 0,$$ it follows
that $${m+n-t-2\choose m-2+i(k-m)+j(l-n)}=0.$$
\end{proof}
\begin{lemma}
Let $a\leq b\leq n$ be positive integers; then $$\dsum_{l=a}^b
{n-l\choose k}={n-a+1\choose k+1}-{n-b\choose k+1}.$$ \label{lem6}
\end{lemma}

\begin{proof}
We prove the equality by induction on $b-a$. If $b-a$, then the
equality is trivially holds. If $b-a\geq 1$, then
$\sum_{l=a}^{b-1} {n-l\choose k}={n-a+1\choose k+1}-{n-b+1\choose
k+1}$ by induction, it follows that $$\ba{rcl} \sum_{l=a}^b
{n-l\choose k} & = & \sum_{l=a}^{b-1} {n-l\choose k}+{n-b\choose
k} \\ & = & {n-a+1\choose k+1}-{n-b+1\choose k+1}+{n-b\choose k}
\\ &  = & {n-a+1\choose k+1}-{n-b\choose k+1} \ea $$
\end{proof}

~~\\ {\bf Proof of Theorem~3.1} We proceed the proof by induction
on $m$. Let $Q_t=(m, t)$ for $t=1, \dots, n-k+1$ (see Figure 2).
Let $S$ be the set of all paths from $P$ to $Q$; then $S$ is the
disjoint union of $S_t$ for $t=1, \dots, n-k$, where $$S_t=\{C\in
S~|~Q_t\in C, Q_{t+1}\notin C\}.$$ For $t=1, \dots, n-k$, let
$Q'_t=(m-1, t)$ and let $S'_t$ be the set of all paths from $P$ to
$Q'_t$. Then there is a one to one correspondence between $S_t$
and $S'_t$ (if $C\in S_t$, then $C-\{Q_t, Q_{t-1}, \dots, Q_1\}\in
S'_t$). Therefore $w(P, Q)=\sum_{t=1}^{n-k} w(P, Q_t')$. To show
the theorem, we consider two situations: (i). $m>l+2$ and (ii).
$m=l+2$. \par Assume first that $m>l+2$. If $t\leq l$, then \be
\label{eq4} w(P, Q'_t)= \sum_{(i, j)\in T} (-1)^{i+j}
{m+n-t-2\choose m-2+i(k-m)+j(l-n)}. \ee by induction. If $t>l$,
then by induction and Remark~\ref{rem2} $$w(P,
Q'_t)={m+n-t-2\choose m-2}-{m+n-t-2\choose k-2}, $$ so that
(\ref{eq4}) also holds in this case by Lemma~\ref{lem5}.
Furthermore, by Lemma~\ref{lem3}(i) and Lemma~\ref{lem4}(i), we
see that \be \label{eq5} \sum_{(i, j)\in T} (-1)^{i+j}
{m+k-2\choose m-1+i(k-m)+j(l-n)}=0.\ee  It follows by
Lemma~\ref{lem6} that $$\ba{rcl}  & & w(P, Q) \\& = &
\sum_{t=1}^{n-k} \sum_{(i, j)\in T} (-1)^{i+j} {m+n-t-2\choose
m-2+i(k-m)+j(l-n)}\\ & = & \sum_{(i, j)\in T} (-1)^{i+j}
\sum_{t=1}^{n-k} {m+n-t-2\choose m-2+i(k-m)+j(l-n)} \\ & = &
\sum_{(i, j)\in T} (-1)^{i+j} [{m+n-2\choose
m-1+i(k-m)+j(l-n)}-{m+k-2\choose m-1+i(k-m)+j(l-n)}] \\ & = &
\sum_{(i, j)\in T} (-1)^{i+j} {m+n-2\choose m-1+i(k-m)+j(l-n)} \ea
$$ Assume now that $m=l+2$. If this is the case, then
$w(P, Q'_1)=w(P, Q'_2)$, so that $w(P, Q)=w(P,
Q'_2)+\sum_{t=2}^{n-k} w(P, Q'_t)$. Again, by induction, $$w(P,
Q'_t)= \sum_{(i, j)\in T} (-1)^{i+j} {m+n-t-2\choose
m-2+i(k-m)+j(l-n)}.$$ Therefore by (\ref{eq5}) and
Lemma~\ref{lem6}, $$\ba{rcl} & &  w(P, Q) \\ & = & \sum_{(i, j)\in
T} (-1)^{i+j} {m+n-4\choose m-2+i(k-m)+j(l-n)}+ \sum_{t=2}^{n-k}
\sum_{(i, j)\in T} (-1)^{i+j} {m+n-t-2\choose m-2+i(k-m)+j(l-n)}
\\ & = & \sum_{(i, j)\in T} (-1)^{i+j}[{m+n-4\choose
m-2+i(k-m)+j(l-n)}+ \sum_{t=2}^{n-k}{m+n-t-2\choose
m-2+i(k-m)+j(l-n)}] \\ & = & \sum_{(i, j)\in T}
(-1)^{i+j}[{m+n-4\choose m-2+i(k-m)+j(l-n)}+ {m+n-3\choose
m-1+i(k-m)+j(l-n)}-{m+k-2\choose m-1+i(k-m)+j(l-n)}] \\ & = &
\sum_{(i, j)\in T} (-1)^{i+j}[{m+n-4\choose m-2+i(k-m)+j(l-n)}
+{m+n-3\choose m-1+i(k-m)+j(l-n)}] \\ & = &  \sum_{(i, j)\in T}
(-1)^{i+j}[{m+n-2\choose m-1+i(k-m)+j(l-n)}-{m+n-4\choose
m-3+i(k-m)+j(l-n)}].  \ea $$ However, by Lemma~\ref{lem3}(ii) and
Lemma~\ref{lem4}(ii) we see that $$\sum_{(i, j)\in T} (-1)^{i+j}
{m+n-4\choose m-3+i(k-m)+j(l-n)}=0.$$ The formula follows. \vskip
0.3in   There are several consequences of Theorem~\ref{theo7}.
\begin{coll}
Let $m$ be a positive integer. Then $$ \sum_{(i, j)\in T}
(-1)^{i+j}{2m\choose m-2i-2j}=2^m,$$ where $T=\{(i, j)\in
\mathbb{Z}\times \mathbb{Z} ~|~i-j=0~or~1\}$.  \label{coll2}
\end{coll}

\begin{proof}
Set $n=m$, $k=m-2$ and $l=n-2$ in the formula of
Theorem~\ref{theo7}, we obtain from Figure 3(a) that $$ \sum_{(i,
j)\in T} (-1)^{i+j}{2m-2\choose m-1-2i-2j}=2^{m-1}.$$ Therefore
the equation holds.
\end{proof}

\begin{picture}(0,0)
\put(0,-10){\line(0, -1){10}} \put(0,-20){\line(1,0){10}}
\put(10,-20){\line(0,-1){10}} \put(10,-30){\line(1,0){10}}

\put(0,-10){\line(1,0){10}} \put(10,-10){\line(0,-1){10}}
\put(10,-20){\line(1,0){10}} \put(20,-20){\line(0,-1){10}}

\put(60,-70){\line(1,0){10}} \put(70,-70){\line(0,-1){10}}
\put(70,-80){\line(1,0){10}} \put(80,-80){\line(0,-1){10}}
\put(60,-70){\line(0,-1){10}} \put(60,-80){\line(1,0){10}}
\put(70,-80){\line(0,-1){10}} \put(70,-90){\line(1,0){10}}
\put(32,-48){$\ddots$}

\put(-3,-13){$\times$} \put(-10,-6){P(1,m)} \put(77,-93){$\times$}
\put(72,-100){Q(m,1)}

\put(20,-120){Figure~3(a)}

\put(130,-10){\line(1,0){20}} \put(150,-10){\line(0,-1){10}}
\put(150,-20){\line(1,0){10}} \put(160,-20){\line(0,-1){10}}
\put(130,-10){\line(0,-1){10}} \put(130,-20){\line(1,0){10}}
\put(140,-20){\line(0,-1){10}} \put(140,-30){\line(1,0){10}}
\put(150,-30){\line(0,-1){10}} \put(160,-60){$\ddots$}
\put(170,-50){$\ddots$}

\put(180,-70){\line(1,0){10}} \put(190,-70){\line(0,-1){10}}
\put(190,-80){\line(1,0){10}} \put(200,-80){\line(0,-1){10}}
\put(200,-90){\line(1,0){10}} \put(210,-90){\line(0,1){20}}
\put(210,-70){\line(-1,0){10}} \put(200,-70){\line(0,1){10}}
\put(200,-60){\line(-1,0){10}}

\put(127,-13){$\times$} \put(120,-7){(1,m)}
\put(207,-93){$\times$} \put(200,-100){(m,1)}

\put(155,-120){Figure~3(b)}

\put(260,-10){\line(1,0){20}} \put(280,-10){\line(0,-1){10}}
\put(280,-20){\line(1,0){10}} \put(290,-20){\line(0,-1){10}}
\put(260,-10){\line(0,-1){10}} \put(260,-20){\line(1,0){10}}
\put(270,-20){\line(0,-1){10}} \put(270,-30){\line(1,0){10}}
\put(280,-30){\line(0,-1){10}} \put(290,-60){$\ddots$}
\put(300,-50){$\ddots$}

\put(310,-70){\line(1,0){10}} \put(320,-70){\line(0,-1){10}}
\put(320,-80){\line(1,0){20}} \put(340,-80){\line(0,1){10}}
\put(340,-70){\line(-1,0){10}} \put(330,-70){\line(0,1){10}}
\put(330,-60){\line(-1,0){10}}

\put(257,-13){$\times$} \put(250,-7){(1,m-1)}
\put(337,-83){$\times$} \put(330,-90){(m,1)}

\put(290,-120){Figure~3(c)}

\end{picture}

\vskip 2.0in
\begin{coll}
Let $F_i$ be the $i$-th Fibonacci number. Then $$F_{2m+1}=
\sum_{(i, j)\in T} (-1)^{i+j}{2m\choose m-3i-2j} $$ and $$F_{2m}=
\sum_{(i, j)\in T} (-1)^{i+j}{2m-1\choose m-3i-2j}, $$ where
$T=\{(i, j)\in \mathbb{Z}\times \mathbb{Z} ~|~i-j=0~or~1\}$.
\label{coll3}
\end{coll}
\begin{proof}
From Figure 3(b), we see that if $n=m$, $k=m-3$ and $l=n-2$ in the
formula of Theorem~\ref{theo7}, then the number of paths from $(1,
m)$ to $(m, 1)$ is  $$\sum_{(i, j)\in T} (-1)^{i+j}{2m-2\choose
m-1-3i-2j}.$$ Moreover, it is easy to see that the path from $(1,
m)$ to $(a, m-a+1)$ is $F_{2a-1}$ for every $a$, therefore
$$F_{2m-1}= \sum_{(i, j)\in T} (-1)^{i+j}{2m-2\choose
m-1-3i-2j}.$$ Similarly, from Figure 3(c), we see that if $n=m-1$,
$k=m-3$ and $l=n-2$ in the formula of Theorem~\ref{theo7}, then
the number of paths from $(1, n)$ to $(m, 1)$ is $F_{2m-2}$, i.e.,
$$F_{2m-2}= \sum_{(i, j)\in T} (-1)^{i+j}{2m-3\choose
m-1-3i-2j}.$$
\end{proof}
\begin{remark}
Corollary~\ref{colll3} gives a new proof to a combinatorial
interpretation of Fibonacci number obtained by Andrews \cite{an}
in 1974.
\end{remark}

\section{Multiplicity of  determinantal ideals}
Let $K$ be a field and $X=(x_{ij})$ be an $m\times n$ matrix of
indeterminates over $K$. Recall that a ladder of $X$ is a subset
$Y$ which satisfies the condition that whenever $x_{ij},
x_{i'j'}\in Y$ with $i<i'$ and $j<j'$, then $x_{ij'}, x_{i'j}\in
Y$. With the help of the previous section, we can find the
multiplicity of certain ladder determinantal ideals.

\begin{theo}
\label{theo8} Let $m, n, l$ and $k$ be non-negative integers with
$k, l\leq min\{m-2, n-2\}$. Let $K$ be a field and $X=(x_{ij})$ be
an $m\times n$ matrix of indeterminats. Let $$Y=\{x_{ij}~|~1\leq
i\leq m, max\{1, l-i+2\}\leq j\leq min\{n, n+m-k-i\}\}$$ be a
ladder of $X$. Let $I=I_{r+1}(Y)$ be the ideal of $K[Y]$ generated
by the $(r+1)$-minors of $Y$ and set $R=K[Y]/I$. Then the
multiplicity of $R$ is $$\det[\sum_{(a, b)\in T}
(-1)^{a+b}{m+n-i-j\choose m-1+(a-b-1)(i-1)+a(k-m)+b(l-n)}
]_{i,j=1, \dots, r},$$ where $T=\{(a, b)\in \mathbb{Z}\times
\mathbb{Z} ~|~a-b=0~or~1\}$.
\end{theo}

\begin{proof}
Let $\tau$ be the lexicographical term order of $K[X]$ induced by
the variable order $$x_{11}>x_{12}> \cdots >x_{1n}>x_{21}> \cdots
>x_{mn}.$$ Let $I^{\ast}$ be the ideal of $K[Y]$
generated by the leading terms of the polynomials of $I$ with
respect to $\tau$ and let $G$ be the set of all $(r+1)$-minors of
$Y$; then it is shown in \cite[Corollary~4.2]{wa} that $G$ is a
Gr\"{o}bner basis for $I$ with respect to $\tau$, so that
$I^{\ast}$ is generated by the leading terms of the polynomials of
$G$. Moreover,, the multiplicity of $R$ coincides with the one of
$K[Y]/I^{\ast}$ and we can associate to $I^{\ast}$ a simplicial
complex $\Delta$ such that $K[\Delta]=K[Y]/I^{\ast}$, where
$K[\Delta]$ is the face ring of $\Delta$. \\ In the following we
identify $X$ with the set  $$\{(i, j)~|~1\leq i\leq m, 1\leq j\leq
n\}.$$ Further we introduce the partial order $(i, j)\leq (i',
j')$ if $i\geq i'$ and $j\leq j'$. As before, if $P, Q\in X$, then
a path from $P$ to $Q$ is a chain in $X$ with endpoints $P$ and
$Q$. Therefore $Y$ is a two-sided diagonal sub-ladder of $X$
determined by $C_k$ and $\tilde{C}_l$. For any subset $Z$ of $Y$,
we use $\delta(Z)$ for the set of points $(i, j)\in Z$ for which
there is no point $(i', j')\in Z$ with $i'<i$ and $j'<j$.\\ Let
$P_i=(i, n)$ and $Q_i=(m, i)$ for $i=1, \dots, r$. Let $F$ be a
facet of $\Delta$; then $F=C_1\cup \cdots \cup C_r$, where
$C_1=\delta (F)$ and $C_i=\delta (F\setminus \cup_{j<i} C_j)$ for
$i\geq 2$. Notice that $C_i\cap C_j=\emptyset$ if $i\neq j$ and
$C_i$ is a path from $P_i$ to $Q_i$ for every $i$. Therefore, $F$
can be decomposed uniquely as a union of non-intersecting paths
from $P_i$ to $Q_i$, $i=1, \dots, r$. Since the union of a family
of non-intersecting paths from $P_i$ to $Q_i$, $i=1, \dots, r$ is
a facet of $\Delta$, we see by \cite[Theorem~5.1.7]{brhe} that the
multiplicity of $K[\Delta]$ is the number of non-intersecting
paths from $P_i$ to $Q_i$, $i=1, \dots, r$. Now, if we replace $m,
n$ and $l$ by $m-i+1$, $n-j+1$ and $l-i-j+2$, respectively, in the
formula of Theorem~\ref{theo7}, then we obtain  $$\ba{rcl}  w(P_i,
Q_j) & = &  \sum_{(a, b)\in T} (-1)^{a+b}{m+n-i-j\choose
m-i+a(k-m+i-1)+b(l-n-i+1)} \\ & = &  \sum_{(a, b)\in T}
(-1)^{a+b}{m+n-i-j\choose m-1+(a-b-1)(i-1)+a(k-m)+b(l-n)}, \ea $$
therefore the multiplicity of $R$ is $$ \det[\sum_{(a, b)\in T}
(-1)^{a+b}{m+n-i-j\choose m-1+(a-b-1)(i-1)+a(k-m)+b(l-n)}]_{i,j=1,
\dots, r}$$ by Theorem~\ref{theo6}.
\end{proof} ~~\\
By Remark~\ref{rem2} we have the following result for one-sided
ladder.

\begin{coll}
Let $k, m$ and $n$ be non-negative integers with $k\leq min\{m-2,
n-2\}$. Let $K$ be a field and $X=(x_{ij})$ be an $m\times n$
matrix of indeterminates. Let $Y=\{x_{ij}~|~1\leq i\leq m,1\leq
j\leq min\{n, n+m-k-i\}\}$ be a ladder of $X$. Let $I=I_{r+1}(Y)$
be the ideal of $K[Y]$ generated by the $(r+1)$-minors of $Y$ and
set $R=K[Y]/I$. Then the multiplicity of $R$ is $$
\det[{m+n-i-j\choose m-i}-{m+n-i-j\choose k-1}]_{i,j=1, \dots,
r}.$$ \label{coll4}
\end{coll}

Let $X$ be a skew-symmetric $n\times n$ matrix of the
indeterminates $x_{ij}$, over a field $K$. Let $R=K[X]/Q_r$, where
$Q_r$ is the ideal of $K[X]$ generated by the set of all
$2r$-pfaffians of $X$. Herzog and Trung \cite{hetr} have given a
formula for the multiplicity of $R$. To end this section, we
generalized their result by giving a formula to the multiplicity
for (diagonal) ladder pfaffian ideals. For the definition and
properties  of pfaffian ideals, the reader is referred to
\cite{brhe}, \cite{hetr} and \cite{ne}.

\begin{theo}
\label{theo9} Let $r, l$ and $n$ be positive integers with $r<l<n$
and $2r\leq n$. Let $X$ be a skew-symmetric matrix of the
indeterminates $x_{ij}$, over a field $K$. Let
$$Y=\{x_{ij}~|~max\{1, i-l\}\leq j\leq min\{n, l+i\}, 1\leq i\leq
n\}$$ be a diagonal ladder of $X$. Let $Q_r$ be the ideal of
$K[Y]$ generated by the set of all $2r$-pfaffians of $Y$ and let
$R=K[Y]/Q_r$; then the multiplicity of $R$ is  $$\det[\sum_{(a,
b)\in T} (-1)^{a+b} {2n-2r-i-j\choose n-r-i+a(-r+i-1)+b(r-l-i)}
]_{i, j=1, \dots, r-1}, $$ where $T=\{(a, b)\in \mathbb{Z}\times
\mathbb{Z} ~|~a-b=0~or~1\}$.
\end{theo}

\begin{proof}
Let $\tau$ be the lexicographical term order of $K[X]$ induced by
the variable order $$x_{1n}> x_{1n-1}> \cdots x_{12}> x_{2n}>
\cdots
> x_{n-2,n-1}>x_{n-1,n}.$$ Let $J_r$ be the set of all
$2r$-pfaffians of $Y$; then it is shown in \cite{hetr} that $J_r$
is a Gr\"{o}bner basis for $Q_r$ with respect to $\tau$. Let
$Q_r^{\ast}$ be the ideal of $K[Y]$ generated by the leading terms
of the polynomials of $J_r$; then the multiplicity of $K[Y]/Q_r$
coincided with the one of $K[Y]/Q_r^{\ast}$ and we can associate
to $Q_r^{\ast}$ a simplicial complex $\Delta$ such that
$K[\Delta]=K[Y]/Q_r^{\ast}$, where $K[\Delta]$ is the face ring of
$\Delta$. \vskip 0.1in ~~ \\ In the following we identify $Y$ with
the set $$\{(i, j)~|~1\leq i\leq j\leq min\{n, l+i\}\}.$$ We
introduce a partial order on $Y$ by $(i, j)\leq (i', j')$ if
$i\leq i'$ and $j\leq j'$. If $Z$ is a subset of $Y$, then we use
$\delta(Z)$ for the set of all points $(i, j)\in Z$ for which
there is no point $(i', j')\in Z$ with $i>i'$ and $j<j'$. \\ Let
$P_i=(i, i+1), Q_i=(n-i, n-i+1)$ for $i=1, \dots, r-1$. Let $F$ be
a facet of $\Delta$; then $F=C_1\cup \cdots \cup C_{r-1}$, where
$C_1=\delta(F)$ and $C_i=\delta(F\setminus \cup_{j<i} C_j)$ for
$i\geq 2$.Notice that $C_i\cap C_j=\emptyset$ if $i\neq j$ and
$C_i$ is a path from $P_i$ to $Q_i$ for every $i$. Therefore, $F$
can be decomposed uniquely as a union of non-intersecting paths
from $P_i$ to $Q_i$, $i=1, \dots, r-1$. Since the union of
non-intersecting paths from $P_i$ to $Q_i$, $i=1, \dots, r-1$ is a
facet of $\Delta$, we see from \cite[Theorem~5.1.7]{brhe} that the
multiplicity of $K[\Delta]$ is the number of non-intersecting
paths from $P_i$ to $Q_i$, $i=1, \dots, r-1$. \\ In order to
obtain our formula, we replace $Y$ by its proper sub-region as
follows. Let $T_i=\{(i, j)~|~i+1\leq j\leq r\}\cup \{(j,
n-i+1)~|n-r+1\leq j\leq n-i\}$ and let $Y'=Y\setminus
\cup_{i=1}^{r-1} T_i$ (see Figure 4 below).

\begin{picture}(0,0)

\put(200,-10){\line(-1,0){40}} \put(160,-10){\line(0,-1){10}}
\put(160,-20){\line(-1,0){10}} \put(150,-20){\line(0,-1){10}}

\put(200,-10){\line(0,-1){50}} \put(200,-60){\line(-1,0){10}}
\put(190,-60){\line(0,-1){10}} \put(190,-70){\line(-1,0){10}}

\put(100,-70){\line(1,0){10}} \put(110,-70){\line(0,1){10}}
\put(110,-60){\line(1,0){10}}

\put(100,-70){\line(0,-1){40}} \put(100,-110){\line(1,0){50}}
\put(150,-110){\line(0,1){10}} \put(150,-100){\line(1,0){10}}
\put(160,-100){\line(0,1){10}}

\put(97,-13){$\times$} \put(90,-5){(1,n)} \put(197,-13){$\times$}
\put(200,-5){$Q_1$(n-r,n)}

\put(157,-13){$\times$} \put(152,-5){(n-l,n)}

\put(197,-53){$\times$} \put(203,-50){$Q_{r-1}$}
\put(197,-63){$\times$} \put(203,-65){(n-r,n-r+1)}

\put(97,-73){$\times$} \put(74,-65){(1,l+1)}
\put(97,-113){$\times$} \put(74,-120){$P_1$(1,r+1)}

\put(137,-113){$\times$} \put(127,-120){$P_{r-1}$}
\put(147,-113){$\times$} \put(150,-120){(r,r+1)}

\put(127,-140){Figure~4}

\end{picture}

\vskip 2.0in

~~\\ Now, observe that if $C$ is a path from $P_i$ to $Q_i$, then
$$\{(i, j)~|~i+1\leq j\leq r+1\}\cup \{(j, n-i+1)~|~n-r\leq j\leq
n-i\}\subseteq C, $$ so that $C\setminus T_i$ is a path of from
$P_i'$ to $Q'_i$, where $P'_i=(i, r+1)$ and $Q'_i=(n-r, n-i+1)$.
Furthermore, if $C'$ is a path from $P'_i$ to $Q'_i$; then $C'\cup
T_i$ is a path from $P_i$ to $Q_i$. Therefore we see that the
multiplicity of $R$ is equal to the number of non-intersecting
paths from $P'_i$ to $Q'_i$ ($i=1, \dots, r-1$). Now, by
reflecting the region $Y$ about $x$-axis and by
Theorem~\ref{theo7} (replace $m, n, k$ and $l$ by $n-r-i+1,
n-r-j+1, n-2r$ and $n-l-i-j+1$, respectively), we obtain $$w(P_i',
Q'_i)= \sum_{(a, b)\in T} (-1)^{a+b} {2n-2r-i-j\choose
n-r-i+a(-r+i-1)+b(r-l-i)}, $$ it follows by Theorem~\ref{theo6}
that the multiplicity of $R$ is $$\det[\sum_{(a, b)\in T}
(-1)^{a+b} {2n-2r-i-j\choose n-r-i+a(-r+i-1)+b(r-l-i)} ]_{i, j=1,
\dots, r-1}. $$
\end{proof}

\begin{coll}
Let $X$ be a skew-symmetric $n\times n$ matrix  of the
indeterminates $x_{ij}$, over a field $K$. Let $r$ be a positive
integer such that $2r\leq n$. Let $Q_r$ be the ideal of $K[X]$
generated by the set of all $2r$-pfaffians of $X$ and let
$R=K[X]/Q_r$; then the multiplicity of $R$ is
$$\det[{2n-2r-i-j\choose n-r-i}-{2n-2r-i-j\choose n-2r-1}]_{i,j=1,
\dots, r-1}.$$
\end{coll}

\vskip 0.2in
\begin{proof}
The formula follows by setting $l=n-1$ in the formula of
Theorem~\ref{theo9} and using the fact that $${2n-2r-i-j\choose
n-r-i+a(-r+i-1)+b(r-l-i)}=0$$ unless $(a, b)=(0, 0)$ or $(1, 0)$.
\end{proof}


\vskip 0.2in


\begin{thebibliography}{1}

\bibitem{ab}
S.S.~Abhyankar,
\newblock {\em Combinatoire des tableaux de {Y}oung, vari\'{e}t\'{e}s
determinantielles
 et calcul de
{H}ilbert (R\'{e}dig\'{e} par {A}. {G}alligo)},
\newblock  Rend. Sem. Mat. Univ. Politec Torino. {\bf 42} (1984), 65--88.

\bibitem{an}
G. Andrews,
\newblock {\em Combinatorial analysis and Fibonacci numbers},
\newblock Fibonacci Quart. {\bf 12} (1974), 141--146.




\bibitem{brhe}
W.~Bruns and J.~Herzog,
\newblock {\em Cohen-{M}acaulay rings},
\newblock Cambridge University Press, 1993.

\bibitem{Ga}
A.~Galligo,
\newblock {\em Computations of some Hilbert functions related with
Schubert calculus},
\newblock  LNIM. Vol. 1124, Springer-Verlag, Berlin/New York, 1985.

\bibitem{gevi}
I.M.~Gessel and G. Viennot,
\newblock {\em Binomial determinants, paths, and hook length formulae},
\newblock Adv. in Math. {\bf 58} (1985), 300--321.

\bibitem{gi}
G.Z.~Giambelli,
\newblock {\em Risoluzione del prolema generale numerativo per gli
spazi plurisecanti di
una curva algebraica},
\newblock Mem. Acad. Sci. Torino (2) {\bf 59} (1909), 433--508.


\bibitem{hetr}
J.~Herzog and N.V. Trung,
\newblock {\em Gr\"{o}bner bases and multiplicity of determinantal
and pfaffian ideals},
\newblock Adv. in Math. {\bf 96} (1992), 1--37.


\bibitem{ne}
E.D. Negri,
\newblock {\em Pfaffian ideals of ladders},
\newblock J. Pure Appl. Algebra {\bf 125} (1998), 141--153.


\bibitem{st}
R. Stanley,
\newblock {\em Combinatorics and commutative algebra},
\newblock Birkh\^{a}user, Basel, 1983.
\newpage

\bibitem{st1}
R. Stanley,
\newblock {\em Enumerative combinatorics I},
\newblock Wadsworth, Pacific Grove, California, 1986.

\bibitem{wa}
H.-J.~Wang,
\newblock {\em Gr\"{o}bner bases and determinantal rings},
\newblock preprint.


\end{thebibliography}
\end{document}